\documentclass[10pt,a4paper]{article}
\usepackage{amssymb, amsthm, amsmath} 
\usepackage{graphicx}
\usepackage{latexsym}
\usepackage{comment}
\usepackage{setspace}

\newtheorem{theorem}{Theorem}[section]
\newtheorem{lemma}[theorem]{Lemma}
\newtheorem{Prop}[theorem]{Proposition}
\newtheorem{Def}[theorem]{Definition}
\newtheorem{Remark}[theorem]{Remark}

\title{Local well-posedness for the Zakharov-Kuznetsov equation in Sobolev spaces}
\author{SATOSHI OSAWA\thanks{Department of Mathematics, Graduate School of Science, Kobe University, 1-1 Rokkodai,
Nada-ku, Kobe-shi, Hyogo,  657-8501, Japan.}}

\begin{document}

\maketitle
\begin{abstract}
The initial value problem for two-dimensional Zakharov-Kuznetsov equation on periodic boundary setting is shown to be locally well-posed in $H^s(\mathbb{R} \times \mathbb{T})$ for $9/10<s<1$. We prove this theorem by using bilinear estimates thinking separetely the first variable and the second variable of space.
\end{abstract}


\section{Introduction}

\ 

\normalsize{W}e consider the initial value problem of the following Zakharov-Kuznetsov equation on the product space $\mathbb{R} \times \mathbb{T}$:
\begin{align}
\left\{ 
\begin{array}{l}
\partial_t u +\partial_x \Delta u + u\partial_x u = 0,\\ \tag{ZK}
u(x,y,0) = u_0 (x,y),
\end{array} \ \ \ \ (x,y) \in \mathbb{R} \times \mathbb{T}, \ \ \ t \in \mathbb{R},
\right.
\end{align}
where $u = u(x,y,t)$ is a real-valued function, $\mathbb{T} = \mathbb{R} / 2 \pi \mathbb{Z}$, and $\Delta = \partial^2_x + \partial^2_y$ is the Laplacian. This equation was introduced by Zakharov and Kuznetsov in \cite{Zakharov}, as a model for the propagetion of ionic-acoustic waves in magnetized plasma. The equation (ZK) is one of the expanded equations of  the Korteweg-de Vries (KdV) equation:
$$\partial_t u + \partial_x^3 u + u \partial_x u = 0.$$

In this paper, we prove a property of well-posedness for an initial data whose norm is bounded in Sobolev spaces. In particular, we observe the behavior of the solution along time development. In this paper, we try to prove a property for the Sobolev space of lower index. 

\ 

The KdV equation is completely integrable, and has infinite number of conserved quantities. Compared to that, the Zakharov-Kuznetsov equation is not completely integrable. However, the Zakharov-Kuznetsov equation has at least two conservation quantities for the its flow: $L^2$  and $H^1$ versions, 
\begin{align}
\left\{
\begin{array}{l}
\displaystyle E[u](t) := \frac{1}{2} \int_{\mathbb{R} \times \mathbb{T}} ({|\nabla u (x,y,t)|}^2 - \frac{1}{3} {u(x,y,t)}^3) dxdy  ,  \\
\displaystyle M[u](t) := \int_{\mathbb{R} \times \mathbb{T}} {u(x,y,t)}^2 dxdy . \nonumber
\end{array}
\right.
\end{align}

\ 

In the non-periodic boundary setting, Gr\"{u}nrock and Herr \cite{Grunrock}, and Molinet and Pilod \cite{Molinet} independently proved the local well-posedness in $H^s(\mathbb{R}^2)$ for $s > 1/2$, by using the Fourier restriction norm method \cite{Bourgain} and Strichartz estimates. Subsequently, Kinoshita \cite{Kinoshita} shows the local well-posedness for $H^s(\mathbb{R}^2)$ using elaborate estimates which are based on region decompositions. By using the $L^2$ and $H^1$ conservation laws, these conclusions specially say that the Zakharov-Kuznetsov equation in $\mathbb{R}^2$ is global well-posed in $L^2$ and $H^1$, respectively. Moreover, similar statements are researched from the perspective of global well-posedness. Shan \cite{Shan} obtained the global well-posedness in $H^s(\mathbb{R}^2)$ for $s>5/7$, and Kinoshita \cite{Kinoshita} obtained the global well-posedness in $L^2(\mathbb{R}^2)$. Shan used the I-method based on $H^1$ conservation laws. Kinoshita used the classical Loomis-Whitney inequality and decompositions according to regions and frequencies. 

On the other hand, in the periodic boundary condition setting, Molinet and Pilod \cite{Molinet} showed the local well-posedness in $H^s(\mathbb{R} \times \mathbb{T})$ for $s \ge 1$ and the global well-posedness in $H^s(\mathbb{R} \times \mathbb{T})$ for $s = 1$. They used nearly the same of strategy in $\mathbb{R}^2$, but generally we cannot use elaborate region decomposition for estimates in $\mathbb{R} \times \mathbb{T}$. Similarily, Linares, Panthee, Robert and Tzvetkov \cite{Linares} obtained the well-poseness result in $H^s(\mathbb{T}^2)$. They showed that the initial value problem for the Zakharov-Kuznetsov equation is locally well-posed in $H^s(\mathbb{T}^2)$ for $s>5/3$. Their proof is based on a variant of  Strichartz estimate. More recently, Kinoshita showed the locally well-posedness in $H^s(\mathbb{T}^2)$, for $s>1$ \cite{Kinoshita2}. 

\ 

In this paper, we will prove the following local well-posedness of initial value problem (ZK) in Sobolev spaces $H^{s}(\mathbb{R} \times \mathbb{T})$ of lower index $s$. 
\begin{theorem}
Let $s > 9/10$. For any $u_0 \in H^s(\mathbb{R} \times \mathbb{T})$, there exists $T=T(\| u_0 \|_{H^s(\mathbb{R} \times \mathbb{T})})$ and a unique solution $u$ of \rm{(ZK)} \it{such that}
$$u \in C([0,T] : H^s((\mathbb{R} \times \mathbb{T})) \cap X^{s, \frac{1}{2}+}_{T}.$$
Moreover, for all $0<T' < T$, there exists a neighborhood $\mathcal{U}$ of $u_0$ in $H^s(\mathbb{R} \times \mathbb{T})$ such that the data-to-solution map 
$$S : v_0 \in \mathcal{U} \mapsto v \in C([0,T'] : H^s((\mathbb{R} \times \mathbb{T})) \cap X^{s, \frac{1}{2}+}_{T'}$$
is smooth, when v is a unique solution of (ZK) to the initial data $v_0$.
\end{theorem}

\ 

In the Section 2, we will show bilinear estimates. We will introduce some lemmas which show relationship between measure of set and regularity of functions, given in \cite{Molinet2}. Moreover, by considering a polynomial called the resonance function, we will confirm strategy of calculation. In Section 3, we will prove Theorem 1.1 by using the propositions discussed in Section 2. We use the frequencies decomposition to observe nonlinear interactions of two wave functions. we calculate the norm of auxiliary space by thinking separately the first variable and the second variable of space. At the higher frequency part, we improved calculation because of Lemma2.9. Finally, we applied Duhamel's principal and evaluated the Sobolev norm of function solution. As a result, we obtain the estimate from above for lower index $s$.

\

\section{Preparation for estimates}
\subsection{Function spaces}

In this subsection, we will introduce some basic function spaces. We aim to prove the initial value problem in Lebesgue spaces or Sobolev spaces. First, we discribe these spaces using the Fourier transform.

\ 

Let $f(x,y)$ be defined on $\mathbb{R} \times \mathbb{T}$. We denote the Fourier transform of $f$ with respect to the variables $(x,y)$ as $\mathcal{F}_{xy}(f) (\xi, q)$. More precisely, $y$ is a variable from $\mathbb{T}$, we should explain $q$ as a variable connecting with the Fourier coefficients. Therefore, we note that $\xi \in \mathbb{R}$ and $q \in \mathbb{Z}$. Moreover, let $\mathcal{F}^{-1}_{xy}$ be the Fourier inverse transform with respect to spacial variables. We use the absolute value of $(\xi,q)$ as $|(\xi,q)|^2 = 3\xi^2 + q^2$ in this paper.  And $a \sim b$ means that there is a constant $c$ such that $a = cb$ and $a \lesssim b$ means that there is a constant $c$ such that $a \le cb$. Moreover, $b+$ in statements means that there exists $\delta >0$ such that $b+\delta$.

Then we define the Sobolev spaces for this discription: equipped with the norm
$$\| f \|^2_{H^s} := \sum_{q \in \mathbb{Z}} \int_{\mathbb{R}} ({\langle |(\xi, q)| \rangle^2 })^s | \mathcal{F}_{xy} (f)(\xi,q)|^2 d\xi ,$$
where  $\langle x \rangle = 1 + |x|$.

\ 

We are able to measure the size and regularity of solutions by using Sobolev spaces. Following \cite{Molinet}, we introduce other more complicated function spaces $X^{s,b}$, called Bourgain's spaces. Now let $u(x,y,t)$ be defined on $\mathbb{R} \times \mathbb{T} \times \mathbb{R}$, and let $\hat{u} (\xi, q, \tau)$ be Fourier tranceform of $u(x, y, t)$ in a similar manner as well. Let $\| u \|_{L^p}$ be the usual Lebesgue norm, and $\| u \|_{l^p}$ be norm in the usual sequence space. \\

\ 

\begin{Def}
$\it{For}$ $s,b \in \mathbb{R}, T>0,$

$$\| u \|_{X^{s,b}} := \| \langle \tau - w(\xi ,q) \rangle^b  \langle |(\xi ,q)| \rangle^{s} |\hat{u}(\xi ,q,\tau)| \|_{L_{\xi ,\tau}^2 l_{q}^{2}},  \eqno{(2.1)} $$ 
$$ \| u \|_{X^{s,b}_T} :=\inf{ \{\| \bar{u} \|_{X^{s,b}} ;\ \  \bar{u} : \mathbb{R} \times \mathbb{T} \rightarrow \mathbb{C},\ \  \bar{u} |_{\mathbb{R} \times \mathbb{T} \times [0,T]} = u \}}, \eqno{(2.2)} $$ 
where $w(\xi,q) = \xi^3 + \xi {q}^2$. When a norm of $u$ introduced \rm{(2.1)} \it{is} bounded, we denote $u \in X^{s,b}$.
\end{Def}
The space $X^{s,b}_T$ will be used for proof of the local well-posedness result, that is under restriction of time on $[0,T]$.
\begin{Remark}
The space $X^{s,b}$ will be characterized by this formula:
$$\| f \|_{X^{s,b}} = \| e^{t \partial_x \Delta} f \|_{H^{b}_{t} H^{s}_{xy}}, \eqno{(2.3)}$$
where $e^{-t \partial_x \Delta}$ is the operator associated linear problem of the Zakharov-Kuznetsov equation, discribed as
$$\mathcal{F}_{xy}((e^{-t \partial_x \Delta} f)) (\xi, q) = e^{itw(\xi,q)} \mathcal{F}_{xy}(f) (\xi, q). \eqno{(2.4)}  $$
\end{Remark}

There are some formulas which discribe relationship about $X^{s,b}$. We introduce some formulas which was stated in \cite{Ginibre}. Let $\eta \in C^{\infty}_0 (\mathbb{R})$ satisfy $0 \le \eta \le 1$, $ \eta_{[-1, 1]} =1$, and $\rm{supp}(\eta) \subset [-2, 2]$. 

\begin{lemma}
\it{Let} $s \in \mathbb{R}$ and $b > 1/2$. Then
$$\| \eta(t) e^{-t \partial_x \Delta} f \|_{X^{s,b}} \lesssim \| f \|_{H^s}, \eqno{(2.5)}$$
for all $f \in H^s$.
\end{lemma}

\begin{lemma}
\it{Let} $s \in \mathbb{R}$. Then 
$$\left \| \eta(t) \int^t_{0} e^{-(t - t') \partial_x \Delta} g(t') dt' \right \|_{X^{s,\frac{1}{2}+}} \lesssim \| g \|_{X^{s , -\frac{1}{2}+}}, \eqno{(2.6)} $$
for all $g \in X^{s, \frac{1}{2}+}.$
\end{lemma}

\begin{lemma}
\it{For} any $T>0$, $s \in \mathbb{R}$ and for all $-1/2 < b' \le b < 1/2$,
$$\| u \|_{X^{s,b'}_T} \lesssim T^{b-b'} \| u \|_{X^{s,b}_T}. \eqno{(2.7)}$$
\end{lemma}

\ 

The polynomial $w(\xi, q)$ which appeared in (2.1) plays an important role in characterizing of solution. Now, we define the resonance function, discribing interaction of solutions:
\begin{align}
\mathcal{H} (\xi_1, \xi_2, q_1, q_2) &= w(\xi_1 + \xi_2, q_1 + q_2) - w(\xi_1, q_1) - w(\xi_2, q_2)  \tag{2.8}\\
 &= 3\xi_1 \xi_2 (\xi_1 + \xi_2) + \xi_2 q^2_1 + \xi_1 q^2_2 + 2(\xi_1 + \xi_2) q_1 q_2. \nonumber
\end{align}

\subsection{Bilinear estimates}

In this subsection, we show the bilinear estimates. First, we define cut-off functions and the dyadic composition. 

\ 

\noindent For $k \in \mathbb{N}$, we denote
 $$\phi(\xi) = \eta(|\xi|) - \eta(2|\xi|),\ \  \phi_{2^k} (\xi, q) = \phi(2^{-k} |(\xi, q)|), \ \ \psi_{2^k}(\xi, q,\tau) = \phi(2^{-k}(\tau - (\xi^3 + \xi q^2))),$$
where $(\xi, \tau) \in \mathbb{R}^2$ and $q \in \mathbb{Z}$.
Here, we prepare notation of interval as $I_N:=\rm{supp} \it{\phi_N} $.
We define the Littlewood-Paley decomposition as
$$P_N(u) := \mathcal{F}^{-1}_{xy} (\phi_N \mathcal{F}_{xy}(u)),\ \ \ Q_L(u) := \mathcal{F}^{-1} (\psi_L \mathcal{F}(u)). \eqno{(2.9)}.$$ 
Now we give the bilinear estimates that are crucial for proving Theorem 1.1. We denote for $a, b \in \mathbb{R}$, $a \wedge b = \min\{a,b\}, a \vee b = \max\{a,b\}$.

\begin{lemma}[Proposition 3.6 in \cite{Molinet}]
Suppose that $N_1, N_2, L_1, L_2$ are dyadic numbers. Then
\begin{align}
&\| (P_{N_1} Q_{L_1} u) (P_{N_2} Q_{L_2} v) \|_{L^2}  \nonumber \\
&\lesssim (N_1 \wedge N_2) (L_1 \wedge L_2)^{\frac{1}{2}}  \| P_{N_1} Q_{L_1} u \|_{L^2}  \|P_{N_2} Q_{L_2} v \|_{L^2}. \tag{2.10}
\end{align}
Moreover, when $4N_2 \le N_1$ or $4N_1 \le N_2$, we obtain
\begin{align}
&\| (P_{N_1} Q_{L_1} u) (P_{N_2} Q_{L_2} v) \|_{L^2}  \nonumber \\
&\lesssim \frac{(N_1 \wedge N_2)^{\frac{1}{2}}}{N_1 \vee N_2} (L_1 \wedge L_2)^{\frac{1}{2}} (L_1 \vee L_2)^{\frac{1}{2}}  \| P_{N_1} Q_{L_1} u \|_{L^2}  \|P_{N_2} Q_{L_2} v \|_{L^2}. \tag{2.11}
\end{align}

\end{lemma}
Molinet and Pilod shows Lemma 2.6, using some basic lemmas introduced in \cite{Molinet2}.

When $f$ is a smooth function or a polynomial of degree 2, Lemma 2.7 and Lemma 2.8 can be written as below.

\begin{lemma} [Lemma 3.8 in \cite{Molinet}]
Let $I$ and $J$ be two intervals on the real line and $f:J \rightarrow \mathbb{R}$ be a smooth function. Then,
$$\rm{mes} \{ \it{x} \in J: f(x) \in I \} \lesssim \frac{|I|}{\rm{inf}_{\it{\xi \in J}}|\it{f}'(\xi)|}, \eqno{(2.12)}  $$
where $\rm{mes} \it{A}$ is the Lebesgue measure of the set A.
\end{lemma}

\begin{lemma} [Lemma 3.9 in \cite{Molinet}]
Let $a \neq 0$, $b$, $c$ be real numbers and $I$ be an interval on the real line. Then, 
$$\rm{mes} \{ \it{x} \in J: ax^{\rm{2}} + \it{bx} + c \in I \} \lesssim \frac{|I|^{\rm{\frac{1}{2}}}}{|a|^{\rm{\frac{1}{2}}}}.  \eqno{(2.13)} $$
\end{lemma}
\noindent Now we denote a new lemma which applied at  the case $N_2/2 \le N_1 \le 2N_2$. We denote by $R_K$ the Fourier multiplier operator in the variable of $x$ as following:
$$R_K f (x) = \int_{\mathbb{R}} \phi (\xi K^{-1}) \mathcal{F}_x f (\xi) e^{ix\xi} d\xi ,$$
where $\mathcal{F}_x$ is the Fourier transform of $f(x)$ with respect to $x$ only.

\begin{lemma}
Let $K$ be a dyadic number. \it{It} holds that
\begin{align}
&\| R_K \left ((P_{N_1} Q_{L_1} u) (P_{N_2} Q_{L_2} v) \right) \|_{L^2}  \nonumber \\
&\lesssim \frac{(N_1 \wedge N_2)^{\frac{1}{2}}}{K^{\frac{1}{4}}} (L_1 \wedge L_2)^{\frac{1}{2}} (L_1 \vee L_2)^{\frac{1}{4}}  \| P_{N_1} Q_{L_1} u \|_{L^2}  \|P_{N_2} Q_{L_2} v \|_{L^2}, \tag{2.14}
\end{align}
where $K$ is a dyadic number satisfying $|\xi| \in I_K$.
\end{lemma}

\noindent \it{Proof} \rm{We} obtain the following estimate from Cauchy-Schwarz inequality and Plancherel's identity:
\begin{align}
\|R_K ( (P_{N_1} Q_{L_1} u) & (P_{N_2} Q_{L_2} v) )\|_{L^2}  \nonumber \\
&= \|R_K ( (\widehat{P_{N_1} Q_{L_1} u}) * (\widehat{P_{N_2} Q_{L_2} v}) )\|_{L^2} \nonumber \\
&= \left \| R_K \left ( \sum_{q_1 \in \mathbb{Z}} \int_{\mathbb{R} \times \mathbb{R}} \widehat{P_{N_1} Q_{L_1} u} (\xi_1, q_1, \tau_1) \cdot \widehat{P_{N_2} Q_{L_2} v}(\xi - \xi_1,q - q_1, \tau - \tau_1) d\xi_1 \tau_1 \right )\right \|_{L^2} \nonumber \\
&\lesssim \sup_{(\xi, q, \tau)} |A^K_{\xi, q, \tau}|^{\frac{1}{2}} \| P_{N_1} Q_{L_1} u \|_{L^2}  \|P_{N_2} Q_{L_2} v \|_{L^2} \nonumber
\end{align}
where
\begin{eqnarray}
A^K_{\xi, q, \tau} = \{ (\xi_1, q_1, \tau_1) \in \mathbb{R} \times \mathbb{Z} \times \mathbb{R} : |\xi| \sim K,  |(\xi_1, q_1)| \in I_{N_1},  |(\xi - \xi_1, q - q_1)| \in I_{N_2},\nonumber \\
 |\tau_1 - w(\xi_1, q_1)| \in I_{L_1}, |\tau - \tau_1 - w(\xi - \xi_1, q - q_1)| \in I_{L_2} \} \nonumber.
\end{eqnarray}
Observe the definition of $A^K_{\xi, q, \tau}$, we get
$$|A^K_{\xi ,q, \tau}| \lesssim (N_1 \wedge N_2)^2 (L_1 \wedge L_2),$$
which implies the estimate in (2.10).\\
On the other hand, the triangle inequality  yields
\begin{align}
|\tau_1 - w(\xi_1, q_1)| + |\tau - \tau_1 - w(\xi - \xi_1, q - q_1)| &\le |\tau - w(\xi_1, q_1)  - w(\xi - \xi_1, q - q_1)| \nonumber \\
&= |\tau - w(\xi, q) - \mathcal{H} (\xi_1, \xi - \xi_1, q_1, q - q_1)|, \nonumber
\end{align}
we get
$$|A^K_{\xi, q, \tau}| \lesssim (L_1 \wedge L_2) |B^K_{\xi, q, \tau}|, \eqno{(2.15)}$$
where
\begin{align}
B^K_{\xi, q, \tau} = \{ (\xi_1, q_1) \in \mathbb{R} \times \mathbb{Z} : |\xi| \sim K, |(\xi_1, q_1)| \in I_{N_1}, |(\xi - \xi_1, q - q_1)| \in I_{N_2}, \nonumber \\
|\tau - w(\xi, q) - \mathcal{H} (\xi_1, \xi - \xi_1, q_1, q - q_1)| \lesssim L_1 \vee L_2 \}. \tag{2.16}
\end{align}
For estimating $B^K_{\xi, q, \tau}$, we calculate the second order differential of $\mathcal{H}$:
$$\left|\frac{\partial^2 \mathcal{H}}{{\partial \xi_1}^2} (\xi_1, \xi - \xi_1, q_1, q-q_1) \right| = 6|\xi| \sim K, \eqno{(2.17)}$$
for $(\xi_1, q_1) \in B^K_{\xi, q, \tau}$. Let $\widetilde{B}_{\xi, q, \tau}(q_1) = \{ \xi_1 \in \mathbb{R}: (\xi_1, q_1) \in B^K_{\xi, q, \tau} \}$ for fixed $q_1$. Then combining Lemma 2.8 and (2.17) we get
$$|\widetilde{B}_{\xi, q, \tau}(q_1)| \lesssim \frac{(L_1 \vee L_2)^{\frac{1}{2}}}{K^{\frac{1}{2}}}, \eqno{(2.18)} $$
for all $q_1 \in \mathbb{Z}$. Finally, it follows from that
$$|B^K_{\xi, q, \tau}| \lesssim \frac{(N_1 \wedge N_2)}{K^{\frac{1}{2}}} (L_1 \vee L_2).$$
Combining with (2.15), we obtain (2.14). $\Box$

\ 

\section{Local well-posedness}
In this section, we prove Theorem 1.1. 
\begin{Prop}
Let $9/10<s<1$. For all $u,v : \mathbb{R} \times \mathbb{T} \times \mathbb{R} \rightarrow \mathbb{C}$ such that $u,v \in X^{s,\frac{1}{2}+}$, we have
$$\| \partial_x (uv) \|_{X^{s,-\frac{1}{2}+}} \lesssim \| u \|_{X^{s,\frac{1}{2}+}} \| v \|_{X^{s,\frac{1}{2}+}} .\eqno{(3.1)} $$
\end{Prop}
Proposition 3.1 instantly allows us to prove Theorem 1.1, which means local well-posedness in Sobolev spaces.\\

\noindent \it{Proof}.  \rm{By} duality, it suffices to prove the following:
$$J := \sum_{q, q_1 \in \mathbb{Z}} \int \Gamma^{\xi_1, q_1, \tau_1}_{\xi, q, \tau} \  \hat{u}(\xi_1, q_1, \tau_1)\  \hat{v}(\xi_2, q_2, \tau_2)\  \hat{w}(\xi, q, \tau)\ d\nu  \lesssim \| u \|_{L^2} \| v \|_{L^2} \| w \|_{L^2}, \eqno{(3.2)}$$
where
$\hat{u}, \hat{v}, \hat{w}$ are non-negative functions,
$$\Gamma^{\xi_1, q_1, \tau_1}_{\xi, q, \tau} = |\xi| \langle |(\xi, q)| \rangle^{s}  \langle |(\xi_1, q_1)| \rangle^{-s} \langle |(\xi_2, q_2)| \rangle^{-s} \langle \sigma \rangle^{-\frac{1}{2}+2\delta} \langle \sigma_1 \rangle^{-\frac{1}{2}-\delta} \langle \sigma_2 \rangle^{-\frac{1}{2}-\delta},$$
$$\xi = \xi_1 + \xi_2, \ \ \ \ \ q = q_1 + q_2, \ \ \ \ \ \tau = \tau_1 + \tau_2,\ \ \ \ \ d\nu = d\xi d\xi_1 d\tau d\tau_1,$$
$$\sigma = \tau - w(\xi, q),\ \ \ \ \ \sigma_1 = \tau_1 - w(\xi_1, q_1),\ \ \ \ \ \sigma_2 =  \tau_2 - w(\xi_2, q_2).$$

Using dyadic decomposition following (2.9), we rewrite $J$ as the following:
$$J = \sum_{N, N_1, N_2} J_{N, N_1, N_2}, \eqno{(3.3)}$$
where
$$J_{N, N_1, N_2}  = \sum_{q,q_1} \int \Gamma^{\xi_1, q_1, \tau_1}_{\xi, q, \tau}\  \widehat{P_{N_1}u}(\xi_1, q_1, \tau_1) \ \widehat{P_{N_2}v}(\xi_2, q_2, \tau_2) \ \widehat{P_N w}(\xi, q, \tau) d\nu .$$

We decompose $J$ to five parts:
$$J_{LL \rightarrow L} := \sum_{N_1 \le 4, N_2 \le 4, N \le 4} J_{N, N_1, N_2},$$
$$J_{LH \rightarrow H} := \sum_{N_2 \le 4, N_1 \le N_2/4} J_{N, N_1, N_2},$$
$$J_{HL \rightarrow H} := \sum_{N_1 \le 4, N_2 \le N_1/4} J_{N, N_1, N_2},$$
$$J_{HH \rightarrow L} := \sum_{N_1 \le 4, N \le N_1/4 \ or \ 4 \le N_2, N \le N_2/4} J_{N, N_1, N_2},$$
$$J_{HH \rightarrow H} := J - (J_{LL \rightarrow L} + J_{HL \rightarrow H} + J_{LH \rightarrow H} + J_{HH \rightarrow L}).$$

We note that $J_{HH \rightarrow H}$ means the  sum of $J_{N, N_1, N_2}$ at the regions of $N_1 \ge 4, N_2 \ge 4, N_2/2 \le N_1 \le 2N_2, N_1/2 \le N \le 2N_1, N_2/2 \le N \le 2N_2$. Estimates of  $J_{LH \rightarrow H}$, $J_{HL \rightarrow H}$ and $J_{HH \rightarrow L}$  are similar to the case of (ZK) on $\mathbb{R}^2$ \cite{Molinet}. In these cases, we use (2.11) directly. 

\ 

\noindent \underline{(i) \it{Estimate of} $J_{LL \rightarrow L}$}. \rm{We} use Plancherel's identity and H\"{o}lder's inequality to obtain 
\[
J_{LL \rightarrow L} \lesssim 
\sum_{N_1 \le 4, N_2 \le 4, N \le 4} \left \|  \left( \frac{\widehat{P_{N_1}u}}{\langle \sigma_1 \rangle^{\frac{1}{2}+}} \right )^{\vee} \right \|_{L^4} 
\left \| \left (\frac{\widehat{P_{N_2}v}}{\langle  \sigma_2 \rangle^{\frac{1}{2}+}} \right )^{\vee} \right \|_{L^4}
\| P_N w \|_{L^2}, \eqno{(3.4)}
\]
where $u^{\vee} (x, y, t)$ denotes the Fourier inverse transform with respect to $(\xi, q, \tau)$. From the Sobolev embedding $H^{1/4+}(\mathbb{R}) \subset L^4(\mathbb{R})$ and Plancherel's identity, we get
\begin{align}
&\left \|  \left( \frac{\widehat{P_{N_1}u}}{\langle \sigma_1 \rangle^{\frac{1}{2}+}} \right )^{\vee} \right \|_{L^4} \left \| \left (\frac{\widehat{P_{N_2}v}}{\langle  \sigma_2 \rangle^{\frac{1}{2}+}} \right )^{\vee} \right \|_{L^4} \| P_N w \|_{L^2} \nonumber \\
&\lesssim  \| \langle \sigma_1 \rangle^{-\frac{1}{4}+} \widehat{P_{N_1} u} \|_{L^2} \cdot  \| \langle \sigma_2 \rangle^{-\frac{1}{4}+} \widehat{P_{N_2} v} \|_{L^2} \cdot \| P_N w \|_{L^2} \nonumber \\
&\lesssim \| P_{N_1} u \|_{L^2} \| P_{N_2} v \|_{L^2} \| P_N w \|_{L^2}, \nonumber
\end{align}
where $N_1 \le 4, N_2 \le 4, N \le 4$. Thus, 
$$J_{LL \rightarrow L} \lesssim \sum_{N_1 \le 4, N_2 \le 4, N \le 4}  \| P_{N_1} u \|_{L^2} \| P_{N_2} v \|_{L^2} \| P_N w \|_{L^2},$$
which yields
$$J_{LL \rightarrow L} \lesssim \| u \|_{L^2} \| v \|_{L^2} \| w \|_{L^2}.$$

\ 

\noindent \underline{(ii) \it{Estimete of} $J_{LH \rightarrow H}$}. \rm{We} use another dyadic decomposition by $Q_L$. We denote
$$J_{N, N_1, N_2} = \sum_{L, L_1, L_2} J^{L, L_1, L_2}_{N, N_1, N_2},$$
where
$$J^{L, L_1, L_2}_{N, N_1, N_2} = \sum_{q,q_1} \int \Gamma^{\xi_1, q_1, \tau_1}_{\xi, q, \tau}\  \widehat{P_{N_1} Q_{L_1} u}(\xi_1, q_1, \tau_1) \ \widehat{P_{N_2} Q_{L_2} v}(\xi_2, q_2, \tau_2)\  \widehat{P_N Q_L w}(\xi, q, \tau) d\nu. \eqno{(3.5)} $$

\noindent Applying the Cauchy-Schwarz inequality, we have that the contribution of this case to $J^{L, L_1, L_2}_{N, N_1, N_2}$ is bounded by
$$C N_2 N^{-s}_1 L^{-\frac{1}{2}+2 \delta} L^{-\frac{1}{2}-\delta}_1 L^{-\frac{1}{2}-\delta}_2 \| (P_{N_1} Q_{L_1} u) (P_{N_2} Q_{L_2} w)\|_{L^2} \| P_N Q_L w \|_{L^2}. \eqno{(3.6)} $$

\noindent By Lemma 2.6, we can show that
$$\| (P_{N_1} Q_{L_1} u) (P_{N_2} Q_{L_2} w)\|_{L^2} \lesssim \frac{N^{\frac{1}{2}}_1}{N_2} L^\frac{1}{2}_1L^\frac{1}{2}_2  \| P_{N_1} Q_{L_1} u \|_{L^2} \| P_{N_2} Q_{L_2} v \|_{L^2}. \eqno{(3.7)}$$
Therefore, combining (3.6) and (3.7), there is some  $\delta>0$ such that for $s > 1/2$, 

\begin{eqnarray}
J_{LH \rightarrow H} & \lesssim&  \sum_{L, L_1, L_2} L^{-\frac{1}{2}+2\delta} L^{-\delta}_1L^{-\delta}_2 \sum_{4 \le N_2, N_1 \le N_2/4} N^{-s + \frac{1}{2}} \| P_{N_1} Q_{L_1} u \|_{L^2} \| P_{N_2} Q_{L_2} v \|_{L^2} \| P_N Q_L w \|_{L^2}  \nonumber \\
& \lesssim & \| u \|_{L^2} \sum_{N \sim N_2} N^{-s+\frac{1}{2}} \| P_{N_2} v \|_{L2} \| P_N w \|_{L^2}   \nonumber \\
& \le & \| u \|_{L^2} \left ( \sum_{N_2} \| P_{N_2} v \|^2_{L^2} \right )^\frac{1}{2}  \left ( \sum_{N} N^{-2s + 1} \| P_{N} w \|^2_{L^2} \right )^\frac{1}{2}  \nonumber \\
& \lesssim & \| u \|_{L^2} \| v \|_{L^2} \| w \|_{L^2}.  \nonumber
\end{eqnarray}

\ 

\noindent \underline{(iii) \it{Estimate of} $J_{HL \rightarrow H}$}. \rm{In} this case, the proof following directly from the proceeding symmetry argument. Indeed, 
$$J_{HL \rightarrow H} \lesssim \| u \|_{L^2} \| v \|_{L^2} \| w \|_{L^2}.$$

\noindent (iv) \underline{\it{Estimate of} $J_{HH \rightarrow L}$}. \rm{By} symmetry, we may assume $4 \le N_1$ and $N \le N_1/4$. Similar to the case of (ii), we have that the contribution of this case to $J^{L, L_1, L_2}_{N, N_1, N_2}$ is bounded by
$$C \frac{N^{s+1}} {N^s_1N^s_2} L^{-\frac{1}{2}+2\delta} L^{-\frac{1}{2}-\delta}_1 L^{-\frac{1}{2}-\delta}_2 \| \widetilde{(P_{N_1} Q_{L_1} u)} (P_{N} Q_{L} w)\|_{L^2} \| P_{N_2} Q_{L_2} v \|_{L^2}, \eqno{(3.8)} $$ 
where $\tilde{f}(\xi, q, \tau) = f(-\xi, -q, -\tau).$ 
We use the interpolation inequality between (2.10) and (2.11). For any $0 < \theta < 1$, 
\begin{align}
&\| \widetilde{(P_{N_1} Q_{L_1} u)} (P_{N} Q_{L} w)\|_{L^2} \nonumber \\
&\lesssim \frac{(N_1 \wedge N)^{\frac{1}{2}(1-\theta) + \theta}}{(N_1 \vee N)^{1-\theta}} (L_1 \wedge L)^{\frac{1}{2}(1-\theta) + \frac{\theta}{2}} (L_1 \vee L)^{\frac{1}{2}(1-\theta)}  \| P_{N_1} Q_{L_1} u \|_{L^2} \| P_N Q_L w \|_{L^2} \nonumber \\
&\lesssim \frac{N^{\frac{1}{2}(1 + \theta)}}{N^{1 - \theta}_1} (L_1 \vee L)^{\frac{1}{2}(1- \theta)} (L_1 \wedge L)^\frac{1}{2} \| P_{N_1} Q_{L_1} u \|_{L^2} \| P_N Q_L w \|_{L^2}. \tag{3.9} 
\end{align}
Then, we substitute (3.9) into (3.8). The contribution of this case to $J^{L, L_1, L_2}_{N, N_1, N_2}$ is bounded by
\begin{eqnarray}
C N^{\frac{1}{2}+\theta} N^{-s+\theta}_1L^{-\delta}_1 L^{-\frac{1}{2}-\delta}_2 L^{2\delta - \frac{\theta}{2}} \| P_{N_1} Q_{L_1} u \|_{L^2} \| P_N Q_L w \|_{L^2} \| P_{N_2} Q_{L_2} v \|_{L^2} \nonumber \\
\lesssim  N^{-s+\frac{1}{2}+2\theta} L^{-\delta}_1 L^{-\frac{1}{2}-\delta}_2 L^{2\delta - \frac{\theta}{2}} \| P_{N_1} Q_{L_1} u \|_{L^2} \| P_N Q_L w \|_{L^2} \| P_{N_2} Q_{L_2} v \|_{L^2}. \nonumber 
\end{eqnarray}
Thus summing up with respect to $N, N_1 N_2, L, L_1, L_2$, we obtain
$$J_{HH \rightarrow L} \lesssim \| u \|_{L^2} \| v \|_{L^2} \| w \|_{L^2},$$
for $0 < \delta < s/4 \ll 1$ and $s > 1/2 + 2\theta$.\\

\noindent \underline{(v) \it{Estimate of} $J_{HH \rightarrow H}$}.  \rm{In} this case, we divide five cases and estimete the resonance function $\mathcal{H}$ in each situation.

\ 

\noindent \underline{(v-1) \it{Case} of $|\xi| \le 100$}. \rm{Let} $k$ be an integer. We denote
$$J_{N, N_1, N_2} =\sum_{k \ge 0} \sum_{L, L_1, L_2} J^{L, L_1, L_2}_{N, N_1, N_2, k}, $$
where
$$J^{L, L_1, L_2}_{N, N_1, N_2, k} = \sum_{q,q_1} \int_{Y_k} \Gamma^{\xi_1, q_1, \tau_1}_{\xi, q, \tau}\  \widehat{P_{N_1} Q_{L_1} u}(\xi_1, q_1, \tau_1) \ \widehat{P_{N_2} Q_{L_2} v}(\xi_2, q_2, \tau_2)\  \widehat{P_N Q_L w}(\xi, q, \tau) d\nu,$$
with 
$$Y_k = \{ (\xi, \xi_1, \tau, \tau_1) \in \mathbb{R}^4 : 2^{-(k+1)} 100 \le |\xi| \le 2^{-k} 100\}.$$ 
We apply the Cauchy-Schwarz inequality, and we have that the contribution of this case to $J^{L, L_1, L_2}_{N, N_1, N_2, k}$ is bounded by
$$ C 2^{-k} N^{-s}_1 L^{-\frac{1}{2}+2\delta} L^{-\frac{1}{2}-\delta}_1 L^{-\frac{1}{2}-\delta}_2 \| (P_{N_1} Q_{L_1} u) (P_{N_2} Q_{L_2} w)\|_{L^2} \| P_N Q_L w \|_{L^2} .\eqno{(3.10)}$$
For  $2^{-(k+1)} 100 \le |\xi| \le 2^{-k} 100$, we obtain
$$\left |\frac {\partial^2 \mathcal{H}} {\partial \xi^2_1} (\xi_1, \xi - \xi_1, q_1, q - q_1) \right | = 6|\xi| \sim 2^{-k} .$$ 
Using Lemma 2.9, we obtain
$$\|R_K ( (P_{N_1} Q_{L_1} u) (P_{N_2} Q_{L_2} v ))\|_{L^2} \lesssim 2^{\frac{k}{4}} N^{\frac{1}{2}}_1 (L_1 \vee L_2)^{\frac{1}{4}} (L_1 \wedge L_2)^{\frac{1}{2}} \| P_{N_1} Q_{L_1} u \|_{L^2} \| P_{N_2} Q_{L_2} v \|_{L^2}. \eqno{(3.11)}$$
Therefore, combining (3.10) and (3.11), we get that (3.10) is bounded by
$$C 2^{-\frac{3}{4}k} N^{s-\frac{1}{2}}_1 L^{-\frac{1}{2} + 2\delta} L^{-\delta}_1 L^{-\delta}_2 \| P_{N_1} Q_{L_1} u \|_{L^2} \| P_{N_2} Q_{L_2} v \|_{L^2} \| P_N Q_L w \|_{L^2} .$$
Summing up with respect to $N, N_1, N_2, L, L_1, L_2, k$, we have that the contribution of this case to $J_{HH \rightarrow H}$ is bounded by 
$$C \| u \|_{L^2} \| v \|_{L^2} \| w \|_{L^2}.$$

\ 

\noindent \underline{(v-2) \it{Case} of $|\xi| \ge 100$, $|\xi_1| \wedge |\xi_2| \le 100$}. \rm{By} symmetry, we may suppose $|\xi_2| \le |\xi_1|$. First, we calculate the resonance function:
$$\frac{\partial \mathcal{H}}{\partial \xi_1}(\xi_1, \xi-\xi_1, q, q-q_1) = 3 \xi (\xi-2\xi_1) +q(q-2q_1). $$
We will consider the case when $|\frac{\partial \mathcal{H}}{\partial \xi_1}| \ge \xi^2$. By using the Cauchy-Schwarz inequality such as (v-1), we have that the contribution of this case to $J^{L, L_1, L_2}_{N, N_1, N_2}$ is bounded by
$$\sum_K C K N^{-s}_1 L^{-\frac{1}{2}+2\delta} L^{-\frac{1}{2}-\delta}_1 L^{-\frac{1}{2}-\delta}_2 \|R_K ( (P_{N_1} Q_{L_1} u) (P_{N_2} Q_{L_2} w))\|_{L^2} \|R_K P_N Q_L w \|_{L^2} .$$
Recall the proof of Lemma 2.9. Indeed of (2.17), we use $\frac{\partial \mathcal{H}}{\partial \xi_1} \gtrsim K$ along which Lemma 2.7 to obtain
$$|\widetilde{B}_{\xi, q, \tau}(q_1)| \lesssim \frac{(L_1 \vee L_2)^{\frac{1}{2}}}{K^2},$$
for (2.18). Then, 
$$\|R_K  ((P_{N_1} Q_{L_1} u) (P_{N_2} Q_{L_2} v ))\|_{L^2} \lesssim K^{-1} N^{\frac{1}{2}}_1 (L_1 \vee L_2)^{\frac{1}{2}} (L_1 \wedge L_2)^{\frac{1}{2}} \| P_{N_1} Q_{L_1} u \|_{L^2} \| P_{N_2} Q_{L_2} v \|_{L^2}.$$
Combining the Cauchy-Schwarz inequality and bilinear estimate, we obtain that the contribution of this case to $J^{L, L_1, L_2}_{N, N_1, N_2}$ is bounded by
\begin{align}
& \sum_{K \lesssim N} K N^{-s}_1 L^{-\frac{1}{2}+2\delta} L^{-\frac{1}{2}-\delta}_1 L^{-\frac{1}{2}-\delta}_2 K^{-1} N^{\frac{1}{2}}_1 (L_1 \vee L_2)^{\frac{1}{2}} (L_1 \wedge L_2)^{\frac{1}{2}} \nonumber \\ 
&\times \| P_{N_1} Q_{L_1} u \|_{L^2} \| P_{N_2} Q_{L_2} v \|_{L^2} \|R_K P_{N} Q_{L} w \|_{L^2} \nonumber \\
&\lesssim \sum_{K \lesssim N} N^{-s + \frac{1}{2}} L^{-\frac{1}{2}+2\delta} L^{-\delta}_1 L^{-\delta}_2 \| P_{N_1} Q_{L_1} u \|_{L^2} \| P_{N_2} Q_{L_2} v \|_{L^2} \|R_K P_{N} Q_{L} w \|_{L^2}. \nonumber
\end{align}
So summing up with respect to $K, N, N_1, N_2, L, L_1, L_2$, the right-hand side is bounded by $C\| u \|_{L^2} \| v \|_{L^2} \| w \|_{L^2}$ for $s > 1/2$.\\

Next, we consider the case when $|\frac{\partial \mathcal{H}}{\partial \xi_1}| \ll \xi^2$. In this case when $|\xi - \xi_1| \gtrsim 1$, the proof is similar to the case (i) and (ii). Without loss of generality, we may suppose that $|\xi_2| = |\xi - \xi_1| = o(1)$, where $o$ is the Landau symbol. In this case, we have
\begin{eqnarray}
\frac{\partial \mathcal{H}}{\partial \xi_1}(\xi_1, \xi - \xi_1, q_1, q-q_1) &=& 3 \xi (\xi-\xi_1) -3\xi \xi_1 + q(q-q_1) -qq_1 \nonumber \\
&=& 6\xi (\xi - \xi_1) -3 \xi^2 + q(q - q_1) - qq_1\nonumber \\
&=& o(\xi) + o(\xi^2) + q(q-q_1) -qq_1. \nonumber
\end{eqnarray}
To satisfy the hypothesis of resonance function $|\frac{\partial \mathcal{H}}{\partial \xi_1}| \ll \xi^2$, $N \sim N_1 \sim N_2$ and $|\xi - \xi_1| = o(1)$, we need $\min\{|q|, |q_1|, |q_2|\} \gtrsim |\xi|$ to set $q(q-q_1) -qq_1= q(q_2 - q_1) \sim \xi^2$.
Next we calculate the resonance function $\mathcal{H}$,
\begin{eqnarray}
\mathcal{H}(\xi_1, \xi - \xi_1, q_1, q - q_1)=3\xi \xi_1 (\xi - \xi_1) + (\xi - \xi_1)(q^2-(q-q_1)^2) + (q-q_1)\xi_1 (q+q_1). \nonumber
\end{eqnarray}
We easily confirm that
$$|q| \sim |q_1| \sim |q_2| \sim |\xi| \sim N,\  3\xi \xi_1 (\xi - \xi_1) \sim o(\xi^2),\  (\xi-\xi_1)(q^2 - (q-q_1)^2) \sim o(\xi^2).$$

In the case of $|\mathcal{H}| \gtrsim \xi^2$, we note that $\max \{ |\sigma|, |\sigma_1|, |\sigma_2| \} \gtrsim \xi^2.$
We suppose $\max\{|\sigma|, |\sigma_1|, |\sigma_2|\} = |\sigma|$, without loss of generality. Then $\langle \sigma \rangle^{-\frac{1}{2}} |\xi| \sim 1$. The bilinear estimate (2.11) along with the argument similar to the case (ii) implies that the contribution of this case to $J_{HH \rightarrow H}$ is bounded by $C\| u \|_{L^2} \| v \|_{L^2} \| w \|_{L^2}$. \\

On the other hand, in the case of $|\mathcal{H}| \ll \xi^2$, it holds that $|(q-q_1) \xi_1 (q + q_1)| \lesssim o(\xi^2)$. We obtain $q + q_1 = 0$. Write $q_1 = -q, q_2 = q-q_1 = 2q$. We rewrite resonance function as $\mathcal{H} = 3(\xi - \xi_1) (\xi \xi_1 - q^2)$ in this case. We regard this formula as a quadratic equation of $\xi_1$ as the following:
$$\xi \xi^2_1- (\xi^2 + q^2)\xi_1 - \xi q^2 + \frac{\mathcal{H}}{3} = 0,$$
which is solved by factoring
$$\xi_1 = \frac{(\xi^2 - q^2) \pm \sqrt{(\xi^2 + q^2)^2 - 4 \xi^2 q^2 + \frac{4}{3}\xi \mathcal{H}}}{2 \xi}. \eqno{(3.12)}$$
Now, the discussion in \cite{Kenig} tells us that it is suffice that the contribution of this case to $J_{HH \rightarrow H}$ is bounded by $C \| u \|_{L^2} \| v\|_{L^2} \| w \|_{L^2}$ when
$$I := \sup_{(\xi, q, \tau)} \frac{|\xi|^2}{\langle \sigma \rangle^{2(1-\frac{1}{2}-\delta)} N^{2s}} \sum_{q_1} \int_{|\xi - \xi_1| \lesssim 1} \frac{d\xi_1}{\langle \sigma - \mathcal{H} \rangle^{2(\frac{1}{2}+ \delta)}} < \infty . \eqno{(3.13)}$$
It suffice to say $I < \infty$. We use the change of the variable (3.12):
$$d\xi_1 = \frac{d\mathcal{H}}{\pm 3 \sqrt {(\xi^2 - q^2)^2 - \frac{3}{4} \xi \mathcal{H}}}. \eqno{(3.14)}$$
For simplicity, we suppose $\xi > 0$. The estimate in Lemma 7.15 of \cite{Linares-Ponce} shows
\begin{align}
\left | \int_{|\xi - \xi_1| \lesssim 1} \frac{d\xi_1}{\langle \sigma - \mathcal{H} \rangle^{2(\frac{1}{2}+ \delta)}} \right | &= \left | \int_{|\xi - \xi_1| \lesssim 1} \frac{d \mathcal{H}}{\langle \sigma - \mathcal{H} \rangle^{1+ 2\delta}\times \left (\pm 3 \sqrt {(\xi^2 - q^2)^2 - \frac{3}{4} \xi \mathcal{H}} \right) } \right |  \nonumber \\
&\lesssim \int^{\infty}_{-\infty} \frac{d \mathcal{H}}{\langle \mathcal{H} \rangle^{1+2\delta} \frac{\sqrt{3 \xi}}{2} \times 3 \sqrt {\frac{3}{4} \frac{(\xi^2 - q^2)^2}{\xi} - \mathcal{H}}} \nonumber \\
&\lesssim \frac{1}{|\xi|^{\frac{1}{2}} + |\xi^2 - q^2|}. \tag{3.15}
\end{align}
Then
\begin{align}
\frac{|\xi|^2}{\langle \sigma \rangle^{2(1-\frac{1}{2}-\delta)} N^{2s}} \sum_{q_1 = -q} \int_{|\xi - \xi_1| \lesssim 1} \frac{d\xi_1}{\langle \sigma - \mathcal{H} \rangle^{2(\frac{1}{2}+ \delta)}}  &\lesssim \frac{|\xi|^2}{\langle \sigma \rangle^{1 - 2\delta} N^{2s}} \frac{1}{|\xi|^{\frac{1}{2}} + |\xi^2 - q^2_1|} \nonumber \\
& \lesssim \langle \xi \rangle^{2- \frac{1}{2} -2s} = \langle \xi \rangle^{\frac{3}{2} - 2s}, \nonumber
\end{align}
which is bounded for $s>3/4$.

\ 

Next, we consider case where $|\xi|, |\xi_1|, |\xi_2| \ge 100$. In this case, following proof in \cite{Molinet} based on the calculation
 \[
	\left| \frac{\partial \mathcal{H}}{\partial \xi_1} (\xi_1, \xi - \xi_1, q_1, q - q_1) \right|
  = \left ||(\xi, q)|^{2} - |(\xi_1, q_1)|^{2} \right |,
\] we devide the region into two subregions such that
\begin{align}
&\rm{(v-3)}\ \ \  \it{i} \in \{\rm{0},1,2 \}, \ \it{j} \in \{ \rm{0},1,2 \}, \ \ \exists (\it{i},j) \ \ \ ||(\xi_{\it{i}}, q_i)|^{\rm{2}} - |(\xi_{\it{j}}, q_j)|^{\rm{2}}| \ge \it{N}^{\rm{6}/5} L^{\rm{6} \delta},  \tag{3.16} \\
&\rm{(v-4)}\ \ \  \forall \it{i},j \in \{ \rm{0},1,2 \},\ \ \  ||(\xi_{\it{i}}, \it{q}_i)|^{\rm{2}} - |(\xi_{\it{j}}, \it{q}_j)|^{\rm{2}}| \le \it{N}^{\rm{6}/5} \it{L}^{\rm{6} \delta}, \tag{3.17}
\end{align}
where we rewrite $\xi = \xi_0$ and $q = q_0$.

\ 

\noindent \underline{\it{Case} of \rm{(v-3)}} Without loss of generality, we suppose that (3.16) holds for $i=0$, and $j=1$. In a similar way to the case (v-2), the contribution of the case to $J^{L, L_1, L_2}_{N, N_1, N_2}$ is bounded by
$$C N N^{-s}_1 L^{-\frac{1}{2}+2\delta} L^{-\frac{1}{2}-\delta}_1 L^{-\frac{1}{2}-\delta}_2 \frac {(N_1 \wedge N_2)^{\frac{1}{2}}}{N^{\frac{3}{5}} L^{3 \delta}} \| P_{N_1} Q_{L_1} u \|_{L^2} \| P_{N} Q_{L} w \|_{L^2} \| P_{N_2} Q_{L_2} v \|_{L^2} . \eqno{(3.18)}$$
Because of $|\xi| \lesssim N_1$, we deduce
$$\lesssim N^{\frac{9}{10} - s}_1 L^{-\frac{1}{2}-\delta} L^{-\frac{1}{2}-\delta}_{1} L^{-\frac{1}{2}-\delta}_{2}  \|  P_{N_1} Q_{L_1} u \|_{L^2}  \|  P_{N_2} Q_{L_2} v \|_{L^2}  \|  P_{N} Q_{L} w \|_{L^2}. \eqno{(3.19)}$$
We sum up the right hand-side with respect to $N, N_1, N_2, L, L_1, L_2$ that the contribution of the case to $J_{HH \rightarrow H}$ is bounded by
$$C \| u\|_{L^2} \| v\|_{L^2} \| w\|_{L^2},$$
for $s > 9/10$.

\ 

To continue this formula in another region, we divide into two cases again.\\
We can rewrite resonance function $\mathcal{H}$ as
\begin{align}
\mathcal{H} (\xi_1, \xi - \xi_1, q_1, q - q_1)  &= 3 \xi \xi_1 \xi_2 +\xi_1 q^2 - \xi q^2_1 -2 \xi_1 q q_1 + 2 \xi q q_1 \nonumber \\
& = -3\xi \xi_1 \xi_2 + (6\xi \xi_1 \xi_2 +\xi_1 q^2 - \xi q^2_1 -2 \xi_1 q q_1 + 2 \xi q q_1) \nonumber \\
& =: -3 \xi \xi_1 \xi_2 + P(\xi, \xi_1, q, q_1). \nonumber \tag{3.20}
\end{align}
The straightforward calculation gives us the estimation of $P(\xi, \xi_1, q, q_1)$:
\begin{align}
|P(\xi, \xi_1, q, q_1)| \le &|\xi_{i_1}| (||(\xi_{i_1}, q_{i_1})|^2 - |(\xi_{i_0}, q_{i_0})|^2| + ||(\xi_{i_1}, q_{i_1})|^2 - |(\xi_{i_2}, q_{i_2})|^2|) \nonumber \\
&+ |\xi_{i_2}| (||(\xi_{i_2}, q_{i_2})|^2 - |(\xi_{i_0}, q_{i_0})|^2| + ||(\xi_{i_2}, q_{i_2})|^2 - |(\xi_{i_1}, q_{i_1})|^2|) \nonumber \\
&\le 4|\xi_{i_1}| N^{\frac{6}{5}}L^{6\delta} \tag{3.21}
\end{align}
where $\{ i_0, i_1, i_2\} = \{ 0,1,2\}$ such that $|\xi_{i_0}| \le |\xi_{i_1}| \le |\xi_{i_2}|$. So we get 
$$|\mathcal{H} (\xi_1, \xi - \xi_1, q_1, q - q_1) | \ge |\xi_{i_1}|(3|\xi_{i_0}||\xi_{i_2}| - 4N^{\frac{6}{5}}L^{6\delta}). \eqno{(3.22)}$$
We rewrite $\xi_{i_0} = \xi_{min}$, $\xi_{i_1} = \xi_{med}$, and $\xi_{i_2} = \xi_{max}$. We consider two cases according to (3.22) : \\
\begin{align}
\rm{(v-4)}\ \ \ |\xi_{\it{max}}||\xi_{min}| \gg \it{N}^{\rm{6}/5} \it{L}^{\rm{6}\delta},\nonumber \\
\rm{(v-5)}\ \ \ |\xi_{\it{max}}||\xi_{min}| \lesssim \it{N}^{\rm{6}/5} \it{L}^{\rm{6}\delta}. \nonumber
\end{align}

\ 

\noindent \underline{\it{Case} of \rm{(v-4)}} From the hypothesis, we note that
\begin{align}
\max \{ |\sigma|, |\sigma_1|, |\sigma_2| \} &\gtrsim |\mathcal{H}(\xi_1, \xi - \xi_1, q_1, q - q_1)| \tag{3.23} \\
&\gtrsim \min\{ |\xi|, |\xi_1|, |\xi_2| \} \cdot (\max \{ |\xi|, |\xi_1|, |\xi_2| \} )^2 \nonumber.
\end{align}
We consider the case when $|\sigma|$ attains maximal of the left-hand side of the above inequality. Using (3.23), we confirm this inequality:
$$|\xi| N^{\frac{1}{2} - s}_1 \langle \sigma \rangle^{-\frac{1}{2} + 3\delta}\  \lesssim \ |\xi_{max}|^{-6\delta} N^{\frac{1}{2} - s}_{1} |\xi_{min}|^{-\frac{1}{2} + 3\delta}\  \lesssim \  1,$$
for $s>1/2$. Moreover, we use calculation of resonance function:
$$\left |\frac {\partial^2 \mathcal{H}} {\partial \xi^2_1} (\xi_1, \xi - \xi_1, q_1, q - q_1) \right | = 6|\xi| \gtrsim 1. \eqno{(3.24)}$$
In a similar way to the case of (v-2), the contribution of this case to $J_{HH \rightarrow H}$ is bounded by
$$C \| u\|_{L^2} \| v\|_{L^2} \| w\|_{L^2},$$
for $s>1/2$.

\ 

\noindent \underline{\it{Case} of \rm{(v-5)}}
Finally, we treat the case of (v-5). Applying the Fourier multiplier operator $R_K, R_{K_1}, R_{K_2}$ to $w, u$ and $v$, respectively, we get the contribution this case to $J^{L, L_1, L_2}_{N, N_1, N_2}$ is bounded by
$$C \sum_{K, K_1, K_2} \frac{K}{ N^{s}_{1}} L^{-\frac{1}{2}+2\delta} L^{-\frac{1}{2}-\delta}_{1} L^{-\frac{1}{2}-\delta}_{2} \|R_K ( (P_{N_1} Q_{L_1} R_{K_1} u) (P_{N_2} Q_{L_2} R_{K_2} v))\|_{L^2} \|  P_{N} Q_{L} R_K w \|_{L^2}. \eqno{(3.25)}$$
Using the hypothesis and, we consider dyadic numbers of (3.25):
$$\frac{K}{N^{s}_{1}} L^{-\frac{1}{2}+2\delta} L^{-\frac{1}{2}-\delta}_{1} L^{-\frac{1}{2}-\delta}_{2} \lesssim K N^{-s} (K^{-\frac{3}{4}}_{min} K^{-\frac{3}{4}}_{max}(N^{\frac{5}{6}} L^{6\delta})^{\frac{3}{4}}) L^{-\frac{1}{2}+2\delta} L^{-\frac{1}{2}-\delta}_{1} L^{-\frac{1}{2}-\delta}_{2}, \eqno{(3.26)}$$
where $K_{max} \sim |\xi_{max}|$ and $K_{min} \sim |\xi_{min}|$. Notice that $(K^{-3/4}_{min} K^{-3/4}_{max}(N^{6/5} L^{6\delta})^{3/4}) \gtrsim 1$ and
$$\left |\frac {\partial^2 \mathcal{H}} {\partial \xi^2_1} (\xi_1, \xi - \xi_1, q_1, q - q_1) \right | = 6|\xi| \sim K. $$
So we use Lemma 2.9 to obtain
\begin{align}
&\|R_K ( (P_{N_1} Q_{L_1} R_{K_1} u) (P_{N_2} Q_{L_2} R_{K_2} v) )\|_{L^2} \tag{3.27} \\
&\lesssim (K_1 \wedge K_2)^{\frac{1}{2}} (L_1 \wedge L_2)^{\frac{1}{2}} K^{-\frac{1}{4}} (L_1 \vee L_2)^{\frac{1}{4}} \| P_{N_1} Q_{L_1} R_{K_1} u \|_{L^2} \| P_{N_2} Q_{L_2} R_{K_2} v \|_{L^2}.\nonumber
\end{align}
Since $K (K_1 \wedge K_2)^{\frac{1}{2}} K^{-\frac{1}{4}} K^{-\frac{3}{4}}_{min} K^{-\frac{3}{4}}_{max} \lesssim K^{-\frac{1}{4}}_{min}$, we deduce that (3.25) is bounded by
\begin{align}
&C \sum _{K, K_1, K_2} K^{-\frac{1}{4}}_{min} L^{-\frac{1}{2}+\frac{11}{2}\delta} L^{-\delta}_{1} L^{-\delta}_2 (N^{-\frac{5}{6}})^{\frac{3}{4}} N^{-s} \ \times  \nonumber \\
&\| P_{N_1} Q_{L_1} R_{K_1} u \|_{L^2} \| P_{N_2} Q_{L_2} R_{K_2} v \|_{L^2} \| P_{N} Q_{L} R_K w \|_{L^2}. \tag{3.28}
\end{align}
Notice that $N^{-s}_1 (N^{6/5})^{3/4} \sim N^{9/10 - s}$. Therefore, when $s>9/10$, we can confirm that (3.28) is bounded by
$$C L^{-\frac{1}{2}+\frac{11}{2}\delta} L^{-\delta}_1 L^{-\delta}_2 N^{\frac{9}{10} - s} \| P_{N_1} Q_{L_1} u\|_{L^2} \| P_{N_2} Q_{L_2} v\|_{L^2} \| P_N Q_L w\|_{L^2}.$$
Summing over $N_1, N_2, N, L_1, L_2, L$, we have that the contribution of the case is bounded by $C \| u\|_{L^2} \| v\|_{L^2} \| w\|_{L^2}$. So we obtain
$$J_{HH \rightarrow H} \lesssim \| u \|_{L^2} \| v \|_{L^2} \| w \|_{L^2},$$
for $s>9/10$, which complete the proof of Proposition 3.1.\ $\Box$

\ 

\noindent \it{Proof of} Theorem 1.1. \rm{Let} $s > 9/10$, $a > 0$, $0 < \delta \ll 1$ and $0<T<1$. We define
$$Z_{a,T} := \{ u \in X^{s,\frac{1}{2}+\delta}_T : \| u \|_{X^{s,\frac{1}{2}+\delta}_T} < a \},$$
and the operator $\Phi_{u_0}(u)$ as follows:
$$\Phi_{u_0} (u) := \eta(t) e^{-t \partial_x \Delta} u_0 + \eta(t) \int^t_0 e^{-(t-t') \partial_x \Delta} \partial_x(u^2(t')) dt'.$$
We shall show that $\Phi_{u_0}$ is a contraction on $Z_{a,T}$, where each $T < T_0$ for some constance $0 < T_0 < 1$.
Using Lemma 2.3, Lemma 2.4, Lemma 2.5 and Proposition 3.1, it holds
\begin{align}
\| \Phi_{u_0}(u) \|_{X^{s,\frac{1}{2}+\delta}_T} &\le  c\| u_0 \|_{H^s} + c\| \eta(t) \partial_x u^2 (\cdot, \cdot, t) \|_{X^{s,-\frac{1}{2}+\frac{\delta}{2}}_T} \nonumber \\
&\le c \| u_0 \|_{H^s} + cT^{\delta} \| \partial_x u^2 (\cdot, \cdot, t) \|_{X^{s, -\frac{1}{2}+\delta}_T} \nonumber \\
&\le c\| u_0 \|_{H^s} + cT^{\delta} \| u \|^2_{X^{s,\frac{1}{2}+\delta}_T} \nonumber \\
&\le c\| u_0 \|_{H^s} + cT^{\delta} a^2. \tag{3.29}
\end{align}
Setting $c\| u_0 \|_{H^s} = a/2$ and
$$cT^{\delta}a < \frac{1}{3},$$
we have $ \| \Phi_{u_0}(u) \|_{X^{s,\frac{1}{2}}_T} < a $, which concludes that $\Phi_{u_0} : Z_{a,T} \rightarrow Z_{a,T}$.\\
Now we consider any $u, u' \in Z_{a,T}$. Using Lemma 2.3, Lemma 2.4, Lemma 2.5, Proposition 3.1 and (3.29), 
\begin{align}
\| \Phi_{u_0} (u) - \Phi_{u_0}(u') \|_{X^{s, \frac{1}{2}+\delta}}  &= \| \eta(t) \int^t_0 e^{-(t-t') \partial_x \Delta} (u^2 - u'^2)(t') dt' \|_{X^{s,\frac{1}{2}+\delta}} \nonumber \\
&\le cT^{\delta} a \| u  + u'\|_{X^{s, \frac{1}{2}+\delta}_T} \| u - u'\|_{X^{s, \frac{1}{2}+\delta}_T} \nonumber \\
&\le 2cT^{\delta} a \| u - u'\|_{X^{s, \frac{1}{2}+\delta}_T}  < \frac{2}{3} \| u - u'\|_{X^{s, \frac{1}{2}+\delta}}, \tag{3.30}
\end{align}
which implies that $\Phi_{u_0}$ is a contraction map from $Z_{a,T}$ into itself. So we obtain a unique solution $u = \Phi_{u_0} (u)$ in $Z_{a,T}$, i.e., 
$$u(t) = \eta(t) e^{-t \partial_x \Delta} + \eta(t) \int^t_0 e^{-(t-t') \partial_x \Delta} \partial_x(u^2(t')) dt'. \eqno{(3.31)}$$
Finally, we will prove $u \in C([0, T] : H^s(\mathbb{R} \times \mathbb{T}))$. Note that 
$$\| u\|_{X^{s,b}} = \| e^{-t \partial_x \Delta} u\|_{H^b_t H^s_{xy}.}$$
By using the Sobolev embedding theorem $H^{1/2 +}(\mathbb{R}) \hookrightarrow C(\mathbb{R})$ and the fact that $e^{-t \partial_x \Delta}$ is lead to continuous in $L^{\infty}([0,T]:H^s)$, we get
$$u \in C([0, T] : H^s(\mathbb{R} \times \mathbb{T})).$$
Here we complete the proof of Theorem 1.1. $\Box$

\section*{Acknowledgement}
The author would like to express the deepest appreciation to Professor Hideo Takaoka for his kind support and reassuring guidance. Without his help, this paper would not have been possible. In addition, the author want to thank Naoya Takahashi and Rika Takiguchi for their contribution to the discussion.

\it{E-mail adress}, \textup{Satoshi Osawa: \texttt{sohsawa@math.kobe-u.ac.jp}}


\begin{thebibliography}{99}
 
 \bibitem{Molinet2}{\sc L. Molinet, J. C. Saut, N. Tzvetkov},
\newblock{\em Global well-posedness for the KP-II equation on the background of a non-localized solution},
\newblock{Ann. Inst. H. Poincar\'{e}, Anal. Non Lin\'{e}aire},
{\bf 28} (2011), 653--676.
 
\bibitem{Molinet}{\sc L. Molinet, D. Pilod},
\newblock{\em Bilinear Strichartz estimates for the Zakharov-Kuznetsov equation and applications},
\newblock{Ann. Inst. H. Poincar\'{e}, Anal. Non Lineair\'{e}},
{\bf 32} (2015), 347--371.

\bibitem{Ginibre}{\sc J. Ginibre},
\newblock{\em Le probl\`{e}me de Cauchy pour des EDP semi-lin\'{e}aires p\'{e}riodiques en variables d\'{e}space},
\newblock{Ast\'{e}risque},
{\bf 237} (1996), 163--187.

\bibitem{Grunrock}{\sc A. Gr\"{u}nrock, S. Herr},
\newblock{\em The Fourier restriction norm method for the Zakharov-Kuznetsov equation},
\newblock {Discrete and Continuous Dynamical Systems},
{\bf 34} (2014), 2061--2068.

\bibitem{Colliander}{\sc J. Colliander, M. Keel,  G. Staffilani, H. Takaoka, and T. Tao},
\newblock{\em Global well-posedness for KdV in Sobolev spaces of negative index.}
\newblock{Electronic Journal of Differential Equations},
{\bf 26} (2001), 1--7.

\bibitem{Bourgain}{\sc J. Bourgain},
\newblock{\em Fourier transform restriction phenomena for certain lattice subsets and applications to nonlinear evolution equations, part 1: Schr\"{o}dinger Equations},
\newblock{Geometric and Functional Analysis},
{\bf 3} (1993), 209--262.

\bibitem{Kenig}{\sc C. E. Kenig, G. Ponce, L. Vega},
\newblock{\em A bilinear estimate with applications to the KdV equation},
\newblock{J. Amer. Math. Soc},
{\bf 9} (1996), 573--603.

\bibitem{Zakharov}{\sc E. A. Zalharov, V.E.Kuznetsov},
\newblock{\em On three dimensional solitons},
\newblock{Sov. Phys. JEPT},
{\bf 39} (1974), 285--286.

\bibitem{Linares-Ponce}{\sc F. Linares, G. Ponce},
\newblock{\em Introduction to nonlinear dispersive equations},
\newblock{Springer}, (2009).

\bibitem{Kinoshita}{\sc S. Kinoshita},
\newblock{\em Global well-posedness for the Cauchy problem of the Zakharov-Kuznetsov equation in 2D},
\newblock{arXiv:1905.01490v1 [math.AP]}.

\bibitem{Kinoshita2}{\sc S. Kinoshita, R. Schippa}, 
\newblock{\em Loomis-Whitney-type inequalities and low regularity well-posedness of the periodic Zakharov-Kuznetsov equation},
\newblock{arXiv:2003.02694v1 [math.AP]}.

\bibitem{Linares}{\sc F. Linares, M. Panthee, T. Robert, and N. Tzvetkov}
\newblock{\em On the periodic Zakharov-Kuznetsov equation},
\newblock{arXiv:1809.02027v1 [math.AP]}.

\bibitem{Shan}{\sc M. Shan}, 
\newblock{\em Global well-posedness and global attractor for two-dimensional Zakharov-Kuznetsov equation},
\newblock{arXiv:1810.02984v1 [math.AP]}.


\end{thebibliography}
\end{document}